\documentclass{amsart}

\usepackage{amsthm,amsfonts,amsmath,amssymb,latexsym,epsfig}
\usepackage{upref,amssymb,eucal,ae}
\usepackage[all,cmtip]{xy}

\newtheorem{theorem}{Theorem}

\newenvironment{remark}{\medskip \refstepcounter{theorem}
\noindent  {\bf Remark \thetheorem}.\rm}{\,}

\newcounter{spthe}

\def\tg{{\tilde g}}

\def\<{\langle}

\def\>{\rangle}

\def\a{\alpha}

\def\tg{\tilde{g}}
\def\mb#1{{\mathbb #1}}

\begin{document}

\title[Canonical isometric embeddings of projective spaces into spheres]
{Canonical isometric embeddings of projective spaces into spheres}
\author{Santiago R. Simanca}
\address{888 S Douglas Road, Apt 121, Coral Gables, FL 33124, U.S.A.}
\email{srsimanca@gmail.com}

\begin{abstract}
We define inductively isometric embeddings of 
$\mb{P}^n(\mb{R})$ and $\mb{P}^n(\mb{C})$ (with their canonical 
metrics conveniently scaled) into the 
standard unit sphere, which present the former as the restriction of the
latter to the set of real points. Our argument parallels 
the telescopic construction of $\mb{P}^\infty(\mb{R})$, $\mb{P}^\infty(\mb{C})$,
and $\mb{S}^\infty$ in that, for each $n$, it extends the previous embedding
to the attaching cell, which after a suitable renormalization makes it possible 
for the result to have image in the unit sphere.  
\end{abstract}


\maketitle

\section{Isometric embeddings into spheres}
We recall that if $(M^n,g)$ is a Riemannian manifold isometrically immersed
into the standard unit sphere $(\mb{S}^N, \tg)$ in Euclidean space
$(\mb{R}^{N+1}, \| \phantom{.} \|^2)$, if 
$\alpha$ and $H$ are the second fundamental form and mean curvature vector
of the embedding, the scalar curvature $s_g$ of $g$ relates to the 
extrinsic quantities as 
\begin{equation}
s_g  = n(n-1) +\tg(H,H)- \tg(\a,\a) \, .
\label{sce}
\end{equation}
A {\it canonical embedding} will be one that is a critical point of the 
functional  
\begin{equation}
\Pi(M)= \int_{M} \tg(\a,\a) \, d\mu \, , \label{sf} 
\end{equation}
under deformations of the embedding \cite{si2}.
 If one such is also minimal, 
(\ref{sce}) implies that the embedding is a critical point of the
total scalar curvature among metrics on $M$ that can be realized by isometric
embeddings into $\mb{S}^N$, and so if $N$ is sufficiently large,
by the Nash isometric embedding theorem \cite{nash}, we conclude that
the metric $g$ on $M$ is Einstein;
conversely, if $(M,g)$ is an Einstein Riemannian manifold 
isometrically embedded into $\mb{S}^N$ as a minimal submanifold, then the 
embedding is canonical.
 
\section{Canonical minimal embeddings of projective spaces}
On a circle centered at the origin, let us consider the map
\begin{equation} \label{m1}
(x_0:x_1) \mapsto \iota_1(x)= (2x_0 x_1,(x_0^2-x_1^2)) \in \mb{R}^2
\, . 
\end{equation}
Antipodal points are sent to the same image, and since the map underlies the 
local degree map $\theta \mapsto 2\theta$, different set of antipodal points 
are mapped to different images. Since 
$$
\| \iota_1(x)\|^2 =(x_0^2 + x_1^2 )^2 \, ,
$$
and the components functions are all homogeneous harmonic polynomials of
degree two, we obtain a 2-to-1 minimal immersion  
$$
\iota_1 : \mb{S}^1 \rightarrow \mb{S}^1 
$$
into the unit circle, of codimension zero.
Hence, with the metric on $\mb{P}^1(\mb{R})$ 
induced by that on $\mb{S}^1$, we obtain a minimal isometric
embedding identification
$$
\mb{P}^1(\mb{R}) \hookrightarrow \mb{S}^1\, , 
$$
between the domain circle $\mb{P}^1(\mb{R})$ of length $\pi$ and the
range circle with its metric of length $2\pi$.
 
We now proceed by induction. Let us assume that we have defined an isometric 
2-to-1 minimal immersion $\iota_{n-1}$ of the sphere $\mb{S}^{n-1}(r_{n-1})$ of
radius $r_{n-1}$  into $\mb{S}^{N_{n-1}} \subset \mb{R}^{N_{n-1}+1}$, which
descends to an isometric embedding of the quotient $\mb{P}^{n-1}(\mb{R})$ with
the induced metric, and is such that the Euclidean 
$\mb{R}^{N_{n-1}+1}$ norm of $\iota_{n-1}(x)$ satisfies  
$$
\| \iota_{n-1} (x)\|^2 = \frac{1}{r^4_{n-1}}(x_0^2+ \cdots + x_{n-1}^2)^2 \, .  
$$
We set $x=(x',x_n)$, where $x'=(x_0, \ldots, x_{n-1})$. Then, if 
\begin{equation} \label{cons}
 r_n^4=\frac{(n+1)(n^2-1)}{n^2}r_{n-1}^4 \, , \quad   
b^2=\frac{1}{(n^2-1)r_{n-1}^4}\, , \quad a^2=2n(n+1)b^2 \, , 
\end{equation}
respectively, we consider the map  
\begin{equation} \label{mn}
\iota_n(x)=\frac{1}{\sqrt{n+1}}
(\iota_{n-1}(x'),ax_nx_0, \cdots , ax_nx_{n-1},b(x_0^2+\cdots +
x_{n-1}^2 - nx_n^2))\in \mb{R}^{N_{n-1}+n+2}\, . 
\end{equation}
Its components are all quadratic harmonic polynomials, and we have that 
$$
\| \iota_{n} (x)\|^2 = \frac{1}{r^4_n}(x_0^2+ \cdots + x_n^2)^2 \, .  
$$


\begin{theorem}
Let $r_n$ and $N_n$ be the sequences 
\begin{equation} \label{ser} 
r_n^4 =\left(\frac{n+1}{2}\right)^2(n-1)! \, , 
\quad N_n = \frac{1}{2}n(n+3)-1 \, ,   
\end{equation} 
respectively.  
Then the map given 
inductively by {\rm (\ref{m1}), (\ref{mn})} above, 
defines an isometric 2-to-1 minimal immersion 
$$ 
\iota_n : \mb{S}^n(r_n) \rightarrow \mb{S}^{N_n} \,   
$$
that maps the fibers of 
$$
\begin{array}{rcc}
& & \mb{P}^n(\mb{R})  \\
& & \uparrow \\
\mb{Z}/2 & \hookrightarrow & \mb{S}^{n}(r_n) 
\end{array}
$$
injectively into the image, and with the Einstein metric
on $\mb{P}^n(\mb{R})$ induced by that of 
the sphere $\mb{S}^n(r_n)$, the map descends to an 
isometric minimal embedding 
$$ 
\iota_n : \mb{P}^n(\mb{R}) \hookrightarrow \mb{S}^{N_n} \subset \mb{R}^{N_n+1}
\, ,
$$
and the diagram of isometric immersions  
\begin{equation} 
\begin{array}{rcccc}
& & \mb{P}^n(\mb{R}) & \stackrel{\iota_n}{\hookrightarrow} & \mb{S}^{N_n}\\
& & \uparrow & \nearrow & \\
\mb{Z}/2 & \hookrightarrow & \mb{S}^{n}(r_n)& &
\end{array}
\end{equation}
commutes. 
\end{theorem}

\begin{remark}
For $n=2$, we have that  
\begin{equation} \label{m2}
(x_0:x_1:x_2) \mapsto \iota_2(x)=\frac{1}{\sqrt{3}}
 (\iota_1(x_0,x_1), 2x_0 x_2, 2x_1 x_2,
(x_0^2+x_1^2- 2x_2^2)/\sqrt{3}) \in \mb{R}^5\, , 
\end{equation}
and we get an isometric minimal embedding 
$\iota_2 : \mb{P}^2(\mb{R}) \hookrightarrow \mb{S}^4$ where the metric $g$ on
$\mb{P}^2(\mb{R})$ is Einstein of scalar curvature $4/3$ induced by the 
metric on  $\mb{S}^2(\sqrt{3/2})$. The second fundamental form $\alpha$ is
such that $\| \alpha \|^2= 2/3$ pointwise. Thus, 
$$
\frac{1}{4\pi}\int_{\mb{P}^2(\mb{R})} s_g d\mu_g= 1 =
\chi (\mb{P}^2(\mb{R})) \, ,
$$
and
$$
\int_{\mb{P}^2(\mb{R})} 
\| \alpha \|^2 d\mu_g = 2\pi =\Pi (\mb{P}^2(\mb{R})) \, ,
$$
respectively. This is the Veronese surface of \cite{cdck}, its
canonical property proven in \cite{rss2} by showing that the Euler-Lagrange 
equations for the functional (\ref{sf}), developed in \cite{si2}
with complete generality, are satisfied for the 
sphere background. The same argument proves that all of the embeddings above, 
in any dimension, are canonical; all of them are given by quadratic harmonic
polynomials inducing eigenfunctions of the Laplacian on projective space
for the first nonzero eigenvalue \cite{dcwa,tak}.  

For $n=3$, $r_3=2^{\frac{3}{4}}$, and the metric on
$\mb{P}^3(\mb{R})$ is Einstein of scalar curvature 
$s_g=\frac{1}{2^{\frac{3}{2}}}6 $, and volume
$\mu_g = \pi^2 (2^{\frac{3}{4}})^3= \frac{1}{2}\omega_3 
(2^{\frac{3}{4}})^3$. Thus, 
$$
\frac{{\displaystyle \int s_g d\mu_g }}{{\displaystyle 
\left( \int d\mu_g\right)^{1/3}}}= 6 \pi^{\frac{4}{3}}= 
\frac{6 \omega_{3}^{\frac{2}{3}}}{ 2^{\frac{2}{3}}} \, ,
$$
the sigma constant of $\mb{P}^3(\mb{R})$ \cite{bn}.
\end{remark}

We regard $\mb{P}^n(\mb{R})$ as a real $n$-dimensional submanifold
of $\mb{P}^n(\mb{C})$, and extend the embeddings above 
so that they fit as the restriction of the canonical embeddings of 
their complex alter egos to the set of real points.

We write a point in $\mb{S}^3 \subset \mb{C}^2$ as
$z=(z_0,z_1)$. The Fubini-Study metric $g$ on $\mb{P}^1(\mb{C})$ is defined to
make of the fibration
\begin{equation} \label{fib}
\begin{array}{rcc}
\mb{S}^1 & \hookrightarrow & \mb{S}^3 \\ & & \downarrow \\ & & 
\mb{P}^1(\mb{C})
\end{array}
\end{equation}
a Riemannian submersion, the action of $\mb{S}^1$ on $\mb{S}^3$ given by
$e^{i \theta}\cdot z=(e^{i\theta}z_0, e^{i\theta}z_1)$. The sectional
curvature of a normalized section is given by 
$K_g (e_1,e_2)= 1 + 3 | \< \pi^* e_1, J\pi^* e_2\>|^2$, where
$\{ \pi^* e_1, J\pi^* e_2 \}$ is a horizontal lift of the section to 
the fiber, and $J$ is the complex structure in $\mb{C}^2$. 

We consider the map
\begin{equation} \label{m1c}
(z_0:z_1) \mapsto \iota_1(z)= (2z_0 \overline{z}_1,
(z_0 \overline{z}_0-z_1\overline{z}_1)) \in \mb{C}\times \mb{R}=\mb{R}^3 \, ,
\end{equation}
All points on an $\mb{S}^1$ orbit are mapped onto the same image, and different
orbits are mapped to different points. Using the Euclidean norm in the
range, we have that  
$$
\| \iota_1(z)\|^2 =(|z_0|^2 + |z_1|^2 )^2 \, ,
$$
and passing to the quotient, we obtain an isometric embedding identification
$$
\iota_1: \mb{P}^1(\mb{C}) \hookrightarrow \mb{S}^2 
$$
between $\mb{P}^1(\mb{C})$ with the Fubini-Study metric of volume $\pi$,
and $\mb{S}^2$ with volume $4\pi$. It is minimal of codimension zero,
and by construction, it restricts to give
the isometric embedding (\ref{m1}) of $\mb{P}^1(\mb{R})$ as the set of totally
geodesic real points of $\mb{P}^{1}(\mb{C})$ embedded into its image, 
$$
\begin{array}{ccc}
\mb{P}^1(\mb{C})  & \hookrightarrow & \mb{S}^2 \\
\cup &  & \cup \phantom{1}\\
\mb{P}^1(\mb{R})  & \hookrightarrow & \mb{S}^1 
\end{array} \, .
$$ 
Notice that the  
composition of the projection (\ref{fib}) and $\iota_1$ is the 
Hopf map $H: \mb{S}^3 \rightarrow \mb{S}^2$ generator of the homotopy
group $\pi_3(\mb{S}^2)$.
Indeed, if 
$z=(z_0,z_1)$ with $z_0=x_0+iy_0$, $z_1=x_1+iy_1$, 
we have 
$$ 
H(z)= \left\{ \begin{array}{ll}
{\displaystyle (2(x_0 x_1+y_0y_1),2(y_0x_1-x_0y_1), 
 x_0^2+y_0^2-x_1^2-y_1^2) } 
 & \text{if $(x_1,y_1)\neq (0,0)$}\, , \vspace{1mm} \\ 
(0,0,1) & \text{if $(x_1,y_1)=(0,0)$} \, .
\end{array}\right.
$$ 

In order to proceed by induction, and since we are to regard
the real embeddings (\ref{mn}) as the restriction of the embeddings we are 
about to define to the set of real points of $\mb{P}^n(\mb{C})$, we begin by
observing that the scales 
defined by the constants in (\ref{cons}) are already fixed. 
So let us assume that we have defined a map  
$$
\iota_{n-1}: \mb{S}^{2(n-1)+1}(r_{n-1})
 \rightarrow \mb{S}^{M_{n-1}} \subset \mb{R}^{M_{n-1}+1} \, ,
$$
that descends to an embedding of the base of the Riemannian submersion 
\begin{equation} \label{fibn}
\begin{array}{rcccc}
\mb{S}^1 & \hookrightarrow & \mb{S}^{2(n-1)+1}(r_{n-1})& & \\ & & 
\downarrow & &\\ & & \mb{P}^{n-1}(\mb{C}) & 
\stackrel{\iota_{n-1}}{\hookrightarrow} & \mb{S}^{M_{n-1}} 
\end{array} 
\end{equation}
into $\mb{S}^{M_{n-1}}$ with the desired properties. If
$z=(z',z_n)$, where $z'=(z_0, \ldots, z_{n-1})\in \mb{C}^{2(n+1)}$, we 
consider the map
\begin{equation} \label{mnc}
\iota_n(z)=\frac{1}{\sqrt{n+1}}
(\iota_{n-1}(z'),a\overline{z}_n z_0, \cdots , a\overline{z}_n z_{n-1},
b(|z_0|^2+\cdots + |z_{n-1}|^2 - n|z_n|^2))\in \mb{R}^{M_{n-1}+2n+2}\, .
\end{equation}
Then we have that
$$
\| \iota_{n} (z)\|^2 = \frac{1}{r^4_{n}}(|z_0|^2+ \cdots + 
|z_n|^2)^2 \, , 
$$
and we obtain a map of spheres   
$$
\iota_n : \mb{S}^{2n+1}(r_n) \rightarrow \mb{S}^{M_n}\subset \mb{R}^{M_{n-1}
+2n+2}\, ,  
$$
$M_n=M_{n-1}+2n+1$, which descends to define an isometric minimal
embedding of the base of the Riemannian submersion   
\begin{equation} \label{fn}
\begin{array}{rcccc}
\mb{S}^1 & \hookrightarrow & \mb{S}^{2n+1}(r_n)& & \\ & &
\downarrow & &\\ & & \mb{P}^n(\mb{C}) &
\stackrel{\iota_n}{\hookrightarrow} & \mb{S}^{M_n}
\end{array}
\end{equation}
into $\mb{S}^{M_n}$. 

\begin{theorem}
Let $M_n$ be the sequence 
$$
M_n=(n+1)^2 -2 \, ,
$$
and $r_n$ and $N_n$ be the sequences in {\rm (}\ref{ser}{\rm )}. 
Then the map 
$$
\iota_n : \mb{S}^{2n+1}(r_n) \rightarrow \mb{S}^{M_n} \subset \mb{R}^{M_n+1} 
$$
given inductively by {\rm (\ref{m1c}), (\ref{mnc})} above, 
maps the fibers of the fibration {\rm (}\ref{fn}{\rm )} injectively into the   
image, and with the scaled Fubini-Study metric on $\mb{P}^n(\mb{C})$ induced
by the metric on the sphere $\mb{S}^{2n+1}(r_n)$, the map descends to an
isometric minimal embedding
$$
\iota_n : \mb{P}^n(\mb{C}) \hookrightarrow \mb{S}^{M_n} \subset \mb{R}^{M_n+1}
\, , 
$$
which restricts to the embedding given inductively by
{\rm (\ref{m1}), (\ref{mn})}, on the 
totally geodesic subset of real points $\mb{P}^n(\mb{R})$, 
and the diagram  
\begin{equation} \label{fnf}
\begin{array}{rcccc}
\mb{S}^1 & \hookrightarrow & \mb{S}^{2n+1}(r_n)& & \\ 
& & \downarrow & \searrow & \\ 
& & \mb{P}^n(\mb{C}) & \stackrel{\iota_n}{\hookrightarrow} & \mb{S}^{M_n}\\
& & \cup & & \cup \\ 
& & \mb{P}^n(\mb{R}) & \stackrel{\iota_n}{\hookrightarrow} & \mb{S}^{N_n}\\
& & \uparrow & \nearrow  & \\ 
\mb{Z}/2 & \hookrightarrow & \mb{S}^{n}(r_n)& &  
\end{array}
\end{equation}
commutes.
\end{theorem}

\noindent {\bf Acknowledgement.} 
We thank Dennis Sullivan for stimulating conversations.

\end{document}